\documentclass[12pt,leqno]{article}   
\usepackage{epsfig, amssymb, amsmath, amscd, subfigure, amsthm,color}
\usepackage[total={13.5cm,9in}]{geometry}

\newtheorem{The}{Theorem}
\newtheorem{Lem}{Lemma}

\newtheorem*{GCT}{Geodesic Characterization Theorem}

\theoremstyle{remark}
\newtheorem{Remark}{Remark}

\theoremstyle{definition}

\newcommand{\lk}{\operatorname{lk}}

\def\c4t4{$C^{\prime \prime} (4) - T (4)$ }
\def\plane{\mathbb{R}^2}
\def\threespace{\mathbb{R}^3}
\def\doublel{\vartriangle_i \! \! (L)}

\title{A New Proof that Alternating Links are Non-Trivial }
\author{Iain Moffatt\thanks{ 
\newline  Department of Applied Mathematics, Charles University, Malostransk\'{e} nam. 25, 118~00 Praha~1, Czech Republic. 
\newline {\em Current address:} Department of Mathematics and Statistics, University of South Alabama, Mobile, AL 36688, USA. 
\newline {\em Email:} {\tt imoffatt@jaguar.usouthal.edu} }}
\date{}

\begin{document}

\maketitle

\begin{abstract}
We use a simple geometric argument and small cancellation properties of  link groups to prove that alternating links are non-trivial. Unlike most other proofs of this result, this proof uses only classic results in topology and combinatorial group theory.
\end{abstract}


\section{Statement of Results}

A link is said to be {\em trivial} if it is the boundary of a set of embedded, disjoint discs (called spanning discs) in the ambient space $S^3$. The triviality of a link is perhaps more intuitive when expressed in terms of link diagrams: a link is trivial if and only if it admits a diagram which contains no crossings. 
In general, it is   hard to decide if a given link is trivial or not.  We can, however, decide if an alternating link is trivial or not by little more than looking at its diagram. A link diagram is said to be {\em alternating} if, as we travel around each component of the link, we pass over  and under strands of the link  alternately. An {\em alternating link} is a link that admits an alternating diagram. After a straightforward normalization, an alternating diagram represents a non-trivial link if and only if it contains crossings. We provide a new proof of this result here.

There are several approaches  in the literature for showing that alternating links are non-trivial. Most of these approaches rely upon the use of a powerful knot invariant: the determinant in \cite{Ba1930,Cr,BB}; the Alexander polynomial in \cite{Cr1959} and \cite{Mu1958}; the Jones polynomial in \cite{Ka} and the Q-polynomial in \cite{Ki1987}, although a purely geometric proof was given in \cite{MT1991}.
These proofs provide different perspectives as to why the result holds.
The argument presented in this paper differs from these proofs as it uses Dehn's lemma and a solution to the word problem for link groups to show in a very direct way that spanning disks for the link cannot exist.
This approach therefore uses only classic topology and combinatorial group theory to provide a direct and intuitive proof for the non-triviality of alternating links.

\medskip

 We prove the following.
\begin{The} \label{th:main}
If $L$ is a link admitting an alternating projection with crossings, then $L$ is non-trivial.
\end{The}

Our method of proof is to  show that the longitudes of the appropriate links are non-trivial by solving the word problem of their link group. It then follows that the links themselves are non-trivial.

If $L$ is an oriented link with components $L_1, \ldots , L_n$ and $N_i$ is a tubular neighbourhood of $L_i$, a {\em meridian} $\mu_i$ of $L_i$ is a non-separating simple closed curve in $\partial N_i$ that bounds a disc in $N_i$ and a longitude $\lambda_i$ is a simple closed curve in $\partial N_i$ that is homologous to $L_i$ in $N_i$ and null-homologous in the exterior  $S^3 - L_i$.

A standard and well known consequence of Dehn's lemma and the loop theorem is that
a link is trivial if and only if all of its longitudes are trivial in the link group.
This reduces Theorem~\ref{th:main} to the problem of solving the word problem for the  longitudes of the link. We do this by using a simple geometric argument to rewrite the longitudes of the link in a certain normal form with respect to the checker-board colouring of a link projection.
Specifically, we write the longitude as a curve which intersects white regions of the checker-board colouring of the link projection before any black regions. This normal form allows us to apply some basic results in small cancellation theory and  solve the word problem for the link groups, concluding that the longitudes are non-trivial.

I would like to thank J.~Crisp for his very interesting comments.

\section{A Normal Form for the Longitudes}

The {\em checker-board colouring} of a link projection is an assignment of a  colour black or white to each of the regions of the projection in such a way that  adjacent regions are assigned   different colours.
\begin{Lem} \label{lem:zigzag}
The $i$-th longitude $\lambda_i$  of a link $L$ is homotopic to a simple closed curve $J \subset S^3 -L$ such that, in terms of the projection,
all intersections of $J$ with white regions of the checker-board colouring occur before any intersections with black regions, with respect to a chosen base point and orientation.
\end{Lem}

The reader may find it helpful to refer to figure~1 while reading the following proof.

\begin{proof}
Begin by fixing a projection $D$ of $L$.
For convenience assign a label $x_1, \ldots , x_{k}$ to each of the regions of $D$. Then, up to homotopy, a based, oriented loop in the link complement can be described by a word in the alphabet $\mathcal{A} = \{x_i, x_i^{-1} | i=1, \ldots , k \}$ by assigning the letter $x_i$ whenever the loop passes downwards through the region $x_i$, and the letter $x_i^{-1}$ whenever the loop passes upwards through the region $x_i$. Notice that the checker-board colouring induces a colour on each letter in $\mathcal{A}$.

We need to choose a representative of the longitude. To do this we define the {\em $i$-th  double} $\doublel$ of a link $L$ to be the curve determined by a parallel copy of the $i$-th component of its projection (so this is the curve determined by the black-board framing).
We define the length of $\doublel$, $|\doublel |$, to be the number of times it ``passes through'' the projection plane.   

 It is a standard fact (for example see \cite{Rolf}) that the longitude $\lambda_i$ can be represented by $ \doublel \cdot m_i^{-\lk(i,i)} $, where $m_i$ represents the $i$-th meridian (we will specify a representative of the meridian shortly) and $\lk(i,i)$ is the self-linking number of $L_i$.
Since there are two choices for the double $\doublel$ (either side of $L_i$), we may choose $m_i$ and $\doublel$ such that the longitude $\lambda_i$ is represented by a word of the form
$w = (l_1^{\pm 1}l_2^{\mp 1}l_3^{\pm 1} \cdots l_{2n}^{\mp 1})\cdot ( l_1^{\pm 1}a^{\mp 1})^k$, where $k=|\lk (i,i)|$ and $a\in \mathcal{A}$. Notice that $w$ alternates in colour and sign.

Without loss of generality, we may assume that $l_1$ is black. We will also assume that $l_1$ appears with a positive exponent, {\em ie.} the first letter of $w$ is $l_1^{+1}$. A similar argument deals with the $l_1^{-1}$ case. We now describe how to deform this representative to the form required in the statement of the lemma. We split the argument into three steps.

\noindent {\em Step 1.} Begin by wrapping $\doublel$ around the component $L_i$. To do this  fix the first and last intersection points, $l_1$ and $l_{2n}^{-1}$, of $\doublel$ and slide the arc  $l_2^{-1} \cdots l_{2n-1}$  underneath $L_i$ as in figure~1(b). Now, moving inwards along this arc from both ends, fix the next two intersection points and slide the rest of the arc over $L_i$ as in figure~1(c). Continue this process of sliding the arc under and over for as long as is possible.  This procedure gives a word of the form
$(\ell_1  \ell_2^{-1}  \cdots \ell_{2n}^{-1})(l_1^{\pm 1}a^{\mp 1})^k$,
such that $\ell_1,  \ldots , \ell_{n}^{-1}$ are coloured black and  $\ell_{n+1} ,  \ldots , \ell_{2n}^{-1}$ are white. Again this is shown in figure~1(c).

\noindent {\em Step 2.} Next  pull the curve $ (l_1^{\pm 1}a^{\mp 1})^k$  along the component $L_i$ so that the white intersections follow the deformed $i$-th double for as long as possible, and the remaining intersection points lie in regions on the opposite side of the curve $L_i$.
Thus, since $2|\lk(i,i)| \leq |\doublel|$, 
we obtain a curve
$ (\ell_1  \ell_2^{-1}  \cdots \ell_{2n}^{-1} )( \ell_{2n} \cdots \ell_{2n-k}^{\pm 1} a_{1}^{\mp 1} \cdots a_k )$, where $ a_{1}^{\mp 1}, \ldots ,a_k  $ are black. This is indicated in figure~1(d).

This representative of the longitude doubles back upon itself and so we can remove some of the white pairs of intersection points, giving the isotopic curve
$ \ell_1  \ell_2^{-1}  \cdots \ell_{2n-k-1}^{\mp 1}   a_{1}^{\mp 1} \cdots a_k $, as in  figure~1(e).

Now if $2k=|\doublel|$ we are done, otherwise we must move on to step three.

\noindent {\em Step 3.}
All that remains is to move the remaining white intersections $\ell_{n+1}, \ldots \ell_{2n-k-1}^{\mp 1}$ to the end of the arc. Clearly this can be done by a sequence of the moves shown in figure~2. These moves are indicated in figures~1~(e)~(f)~and~(g).

\end{proof}

\begin{Remark}
Stopping after step 2 in the proof shows that the longitude $\lambda_i$ is conjugate to a curve with the properties of $J$ in the lemma. In actual fact this is enough to prove the main theorem, however we prefer the stronger form of the lemma.
\end{Remark}

\begin{figure}[htbp] \label{fig:norm}

\begin{center}
\subfigure[]{\epsfig{file=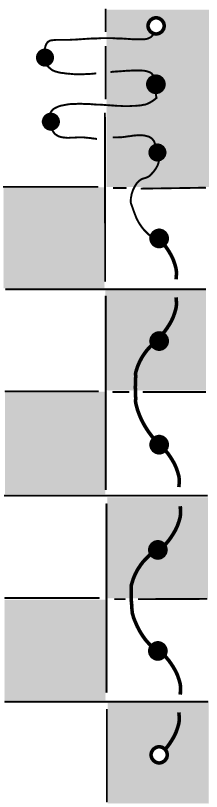, height=5cm}}
\quad \raisebox{2cm}{$\longrightarrow$}  \quad
\subfigure[]{\epsfig{file=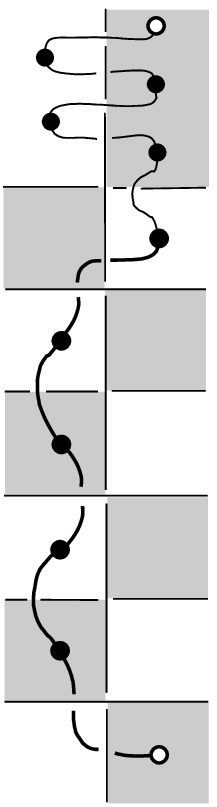, height=5cm}}
\quad \raisebox{2cm}{$\longrightarrow$}  \quad
\subfigure[]{\epsfig{file=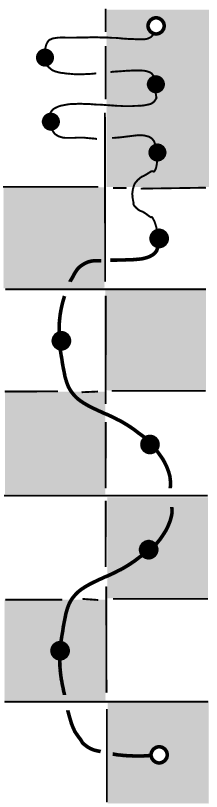, height=5cm}}
\quad \raisebox{2cm}{$\longrightarrow$}
\end{center}

\begin{center}
\subfigure[]{\epsfig{file=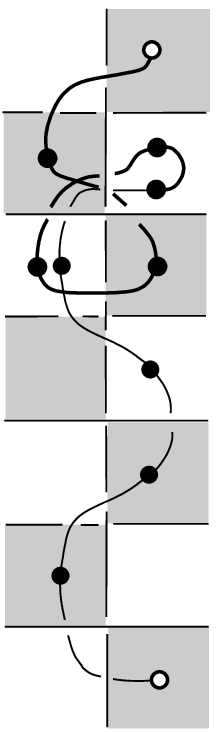, height=5cm}}
\quad \raisebox{2cm}{$\longrightarrow$}  \quad
\subfigure[]{\epsfig{file=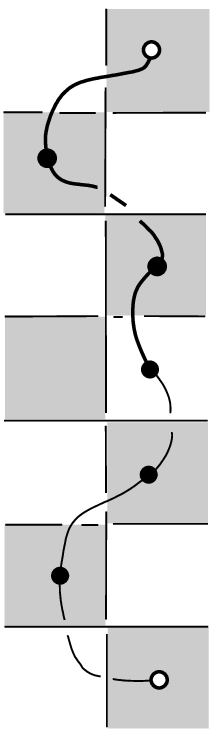, height=5cm}}
\quad  \raisebox{2cm}{$\longrightarrow$} \quad
\subfigure[]{\epsfig{file=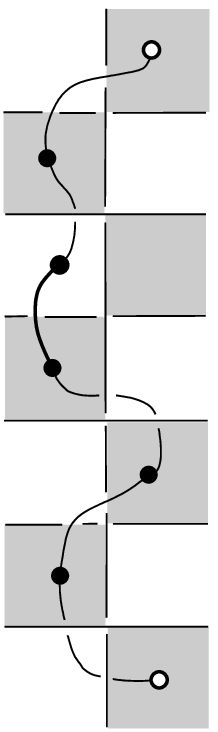, height=5cm}}
\quad \raisebox{2cm}{$\longrightarrow$}  \quad
\subfigure[]{\epsfig{file=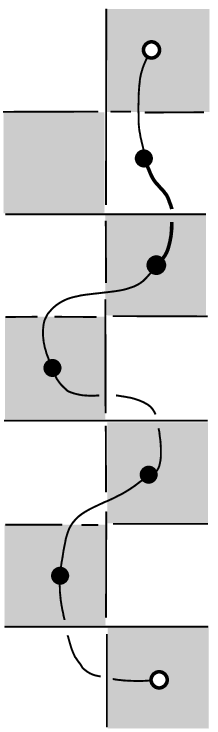, height=5cm}}
\end{center}

\caption{}
\end{figure}

\begin{figure}
\epsfig{file=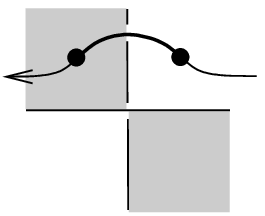, height=1.5cm}
\quad \raisebox{7.5mm}{$\longrightarrow$}  \quad
\epsfig{file=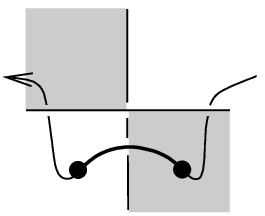, height=1.5cm}
\quad \quad
\raisebox{7.5mm}{and}
\quad \quad
\epsfig{file=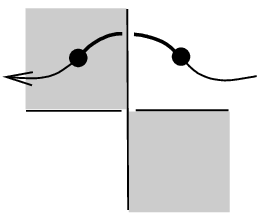, height=1.5cm}
\quad \raisebox{7.5mm}{$\longrightarrow$}  \quad
\epsfig{file=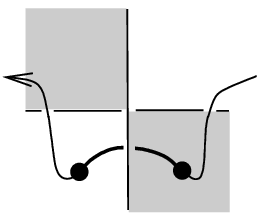, height=1.5cm}
\label{fig:move}
\caption{}
\end{figure}

\section{The Proof of the Theorem}
As in the proof of the lemma, label the regions of the projection $D$ of $L$ with $x_1, \ldots , x_n$.
Now form a group presentation by taking the set of labels $\{ x_1, \ldots ,x_n \}$ as the set of generators and deriving a  relator $x_a x_b^{-1} x_c x_d^{-1}$ from each crossing according to the scheme shown in figure~3.
This gives a presentation which,  after {\em symmetrization}, {\em ie.} adjoining all cyclic permutations and inverses  of the relators to the presentation, we call the {\em augmented Dehn presentation} of $L$. We will see that, although it is an abuse of notation, using the same letters for the regions of the projection and the group generators is natural and convenient.

\begin{figure}
\[\epsfig{file=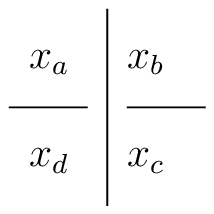, height=2cm}\]
\caption{}
\label{relators}
\end{figure}

We also define the {\em augmented link} obtained from $D$ to be the link obtained by adding an extra unknot component bounding the projection, and forming a link by regarding $\plane$ as the $x$-$y$ plane of $\threespace \cup \infty = S^3$ and  ``pulling the overcrossings up a little.''

If we choose a base point above the plane, a generator $x_i$ of the augmented Dehn presentation is realized in the complement of the augmented link by a loop which passes downwards through the region $x_i \subset \plane \subset \threespace \cup \infty = S^3$ and passes back up through the unbounded region of the  augmented link to the base point.
Notice that up to homotopy there is a clear correspondence between words in the augmented Dehn presentation and words arising from based, oriented loops as in the proof of lemma~\ref{lem:zigzag}.

Clearly the augmented Dehn presentation  is the symmetrization of the Dehn presentation (for example see \cite{LSbook}) of the group of the augmented link. Consequently the augmented Dehn presentation
is a presentation of the free product of the infinite cyclic group and the link group. Therefore solving the word problem for the augmented Dehn presentation of $L$ solves it for the link group of $L$.

It turns out that particular presentations of the groups of certain links have a very strong combinatorial structure. To state this result precisely we need to introduce a little more notation.
A link projection $D$ divides the plane into {\em regions}.
 $D$ is said to be {\em reduced} if at each crossing four distinct regions of the plane meet. Every link admits reduced diagrams. Thus reduced diagrams do not contain either of the configurations
 \[\epsfig{file=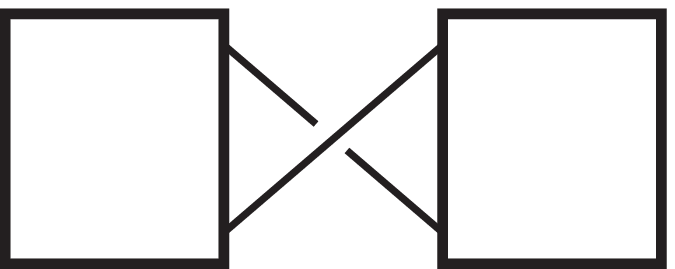, height=1cm} \quad \raisebox{3mm}{\text{ or }} \quad
 \epsfig{file=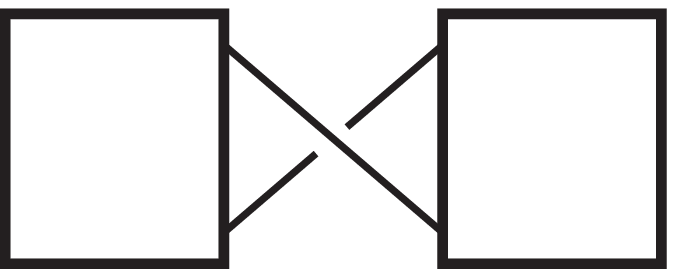, height=1cm},
 \]
 where the rest of the link is contained within the two boxes.
 It is clear from the figures that every diagram is equivalent to a reduced diagram. 
 
A reduced projection is said to be {\em prime} if  it is connected, contains at least one crossing and there does not exist a simple closed curve in the plane intersecting $D$ transversally in exactly two points on different arcs of $D$. A link admits a non-prime diagram if and only if it cannot be expressed as the connected sum of non-trivial links.

Weinbaum proved the following lemma for knots, but his proof also works for links.
\begin{Lem}[Weinbaum~\cite{We1971}]
The augmented Dehn presentation read from  a reduced, prime, alternating projection of a link is a \c4t4 small cancellation group.
\end{Lem}

As we are only interested in the characterization of geodesics, we exclude the definition of a small cancellation group. This can be found in \cite{LSbook}.

A {\em chain} is a Van Kampen diagram having the form shown in figure~4, where $n \geq 1$. We call the word $t_0 t_1 t_2 \cdots t_{n+1}$ a {\em chain word}. (For details on Van Kampen diagrams see \cite{Jo} or \cite{LSbook}. (In \cite{LSbook} they are called ``diagrams."))

\begin{figure}
\[\epsfig{file=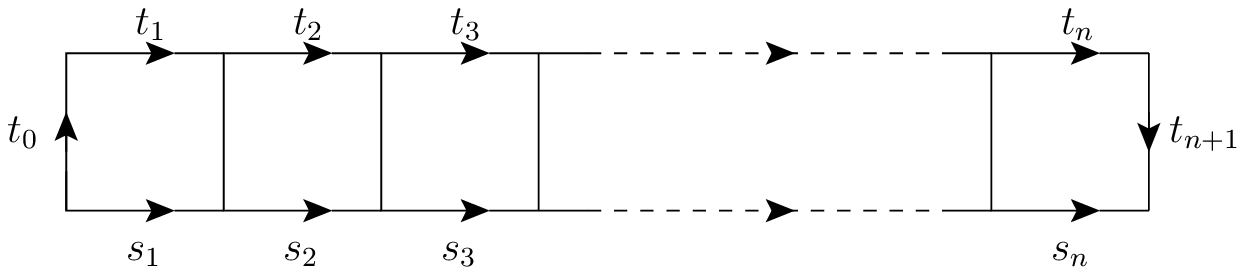, width=10cm}\]
\caption{}
\label{chain}
\end{figure}

Given an arbitrary finite group presentation, the set of lengths of all words representing an element of the group has a minimum. We call any word which attains this minimum a {\em geodesic}. Geodesics in a \c4t4 presentation are characterized by the absence of chain words. For the following theorem, recall that a word is {\em freely reduced} if it contains no subwords of the form $x^{\pm1}x^{\mp 1}$.
\begin{GCT}
A word in a \c4t4 presentation is geodesic if and only if it is freely reduced and contains no chain subwords.
\end{GCT} 

The Geodesic Characterization Theorem appeared implicitly in  \cite{AS1983}  and \cite{GS1990}. A formal proof can be found in \cite{Ka1997}, where the Geodesic Characterization Theorem appears as Lemma~3.2.

Using the checker-board colouring we can assign a {\em parity}, black or white, to each generator-inverse pair according to the colour of the region that generator corresponds to.
Notice that the relators of the augmented Dehn presentation are words which alternate in parity and therefore the horizontal and vertical edges of a chain correspond to letters of different parities.

\smallskip

Putting all this together, we can prove our main result.
\begin{proof}[Proof of Theorem \ref{th:main}]
Since the sum of two non-trivial links is non-trivial (this follows from the additivity of genus under the connect sum operation, for example), it is enough to prove the theorem for prime links.
In this case,  by Lemma~\ref{lem:zigzag} and the geometric interpretation of the generators of the augmented Dehn presentation, the longitude can be represented by  a non-empty word $w$ which changes parity exactly once.
Since the projection is reduced, $w$ is freely reduced.
 A word of this form cannot contain a chain word (as these change parity twice) and since the augmented Dehn presentation is a \c4t4 small cancellation group, the geodesic characterization theorem tells us that the longitudes are non-trivial and therefore the link itself is non-trivial.
\end{proof}

\begin{Remark}
John Crisp has observed that the theory of CAT(0) groups can be used in place of small cancellation theory in proving that the elements of the link group described by Lemma~\ref{lem:zigzag} are non-trivial.

It is also possible to use Dugopolski's algorithm from \cite{Du} to prove that the elements of the link group described by Lemma~\ref{lem:zigzag}. Dugopolski's algorithm  uses normal surface theory rather than group theory.
\end{Remark}

\end{document}